\documentclass{amsart}
\usepackage{amssymb,latexsym}
\usepackage{graphicx}
\input epsf
\newtheorem{theorem}{Theorem}[section]
\newtheorem{proposition}[theorem]{Proposition}
\newtheorem{corollary}[theorem]{Corollary}
\newtheorem{definition}[theorem]{Definition}
\newtheorem{remark}[theorem]{Remark}
\newtheorem{lemma}[theorem]{Lemma}
\newtheorem{example}[theorem]{Example}

\newcommand{\ZZ}{{\mathbb Z}}

\begin{document}
\title [Polygon dissections and  generalized  cluster complexes]
{Polygon dissections and some generalizations of cluster complexes}

\author{Eleni Tzanaki}        
\address{Department of Mathematics\\
University of Crete\\
71409 Heraklion, Crete, Greece}
\email{etzanaki@math.uoc.gr}
\date{\today}
\thanks{The present research will be part of the author's Ph.D. thesis at the University of Crete}

%
\begin{abstract}
   Let $W$ be a Weyl group corresponding to the root system $A_{n-1}$ or $B_n$. 
   We  define a simplicial complex $ \Delta^m_W  $  in terms of polygon dissections
       for such a group  and any  positive integer $m$. 
 For $ m=1 $, $ \Delta^m_W$  is isomorphic to the cluster  complex  corresponding to $ W $, defined  in    \cite{FZ}.
We enumerate the faces of $ \Delta^m_W $  and show that the entries of its  $h$-vector are given by the  generalized Narayana numbers
    $ N^m_W(i)  $, defined in \cite{Atha3}.
      We also prove that for any $ m \geq 1$ the complex $ \Delta^m_W  $ is shellable and hence Cohen-Macaulay.
    \end{abstract}

\maketitle
  \section{Introduction and results}
 \label{intro}
 
   Part of the  motivation for the present work  comes  from the paper \cite{Pr} by J. Przytycki and A. Sikora, 
 where  the number of  certain  polygon dissections 
    was  computed. More specifically, consider the set of dissections of a convex $(mn+2)$-gon $P$ by noncrossing 
  diagonals into polygons each having number of vertices congruent to $2$ modulo $m$. We call the above kind of 
   dissections  \emph{  m-divisible}.

\begin{theorem}\cite[Corollary 2]{Pr}
  \label{note}
   The number of $m$-divisible polygon dissections of an $ (mn+2)$-gon with i diagonals is  equal to 
   \[   \frac{1}{n} \left( \begin{array}{c} mn+i+1 \\ i \end{array} \right)
  \left( \begin{array}{c} n \\ i+1 \end{array} \right), \;\; i=0,\ldots,n-1.\]
\end{theorem}

 Consider the  abstract simplicial complex  $ \Delta^m_{ A_{n-1}}$  whose faces are the $m$-divisible dissections of  an $( mn+2)$-gon $P$. 
Precisely, the vertices of $ \Delta^m_{A_{n-1}}$ are the diagonals which dissect $P$  into an $(mj+2)$-gon and an $(m(n-j)+2)$-gon $(1\leq j \leq n-1)$ 
and the faces are sets of diagonals of this type which are pairwise noncrossing. 
 The facets of   $\Delta^m_{A_{n-1}}$ are the dissections into $ (m+2)$-gons. 
 Since every such dissection   has $ n-1 $ diagonals,  $\Delta^m_{A_{n-1}} $ is  pure of dimension $ n-2 $. 
  For $ m=1$ this  simplicial complex is isomorphic to the boundary complex of the  $(n-1)$-dimensional simplicial associahedron \cite[p. 18]{Zi}. 
   In this language, Theorem \ref{note} computes the number of $( i-1)$-dimensional faces of   $ \Delta^m_{A_{n-1}}$. 

    S. Fomin and  A. Zelevinsky  \cite{FZ}  have introduced a  simplicial complex associated to a  crystallographic root system  $\Phi$,  
  called the cluster complex. By abuse of notation we denote this complex by $\Delta_W$, where  $W$ is  the Weyl group associated to $\Phi $. It was proved in
   \cite{CFZ}  that $\Delta_W$  is isomorphic to the boundary complex  of a  simplicial convex polytope which is  called the simplicial   generalized 
 associahedron  associated to $ \Phi $. 
For $ W=A_{n-1} $ this  polytope is the classical simplicial associahedron. Thus $\Delta_{A_{n-1}}$ 
    is combinatorially isomorphic to $\Delta^m_{A_{n-1}}$ for $m=1$. 
     We define  a complex  $ \Delta^m_{B_n} $  in terms of polygon dissections so that $\Delta_{B_n} $ is 
     isomorphic to $ \Delta^m_{B_n}  $ for $m=1$.   

   The complex $ \Delta^m_{B_n}$ is defined as follows. 
  Let  $ P $ be  a centrally symmetric convex polygon with $ 2mn+2$ vertices labeled  in anticlockwise order
   $1,$ $2,\ldots,$ $mn+1, $ $ {\overline 1},$ ${\overline 2},\ldots,$ $\overline {mn+1}$.  A $B$-diagonal in $ P $ is either   $ \rm(i) $ a diameter, i.e. a diagonal 
   joining antipodal points $ i $ and $ \overline i $ for some  $1\leq i \leq mn+1$, or $ \rm(ii) $ a pair of noncrossing diagonals $ ij, \overline i \overline j $ for two distinct
    $ i,j \in \{1,2,\ldots,mn+1,\overline 1,\overline 2,\ldots,\overline {mn+1}  \} $, nonconsecutive around the boundary 
  of the polygon. (It is understood that if $a={\bar i}$ for some positive integer $i$, then ${\bar a}=i$.)
We define $ \Delta^m_{B_n} $  as  the simplicial complex whose vertices are the $B$-diagonals 
  that dissect $P$ into a pair of $(mj+2)$-gons and a centrally symmetric $( 2m(n-j)+2)$-gon $ (1\leq j \leq n)$
and whose faces are the sets of $B$-diagonals of this type which are pairwise noncrossing. 
The facets of $\Delta^m_{B_n} $ are the dissections of $P$ by $B$-diagonals   into   $(m+2)$-gons. 
 Any  such dissection has $n$ $B$-diagonals and hence  the complex $\Delta^m_{B_n} $ is pure of dimension $ n-1$. 
    For $ m=1 $, $\Delta^m_{B_n} $ is  isomorphic  to  the boundary complex of Simion's type-B associahedron, 
       also known as  cyclohedron \cite{Si} or   Bott-Taubes polytope  \cite{BT}.
       In Section \ref{enumeration} we compute the $f$-vector of $ \Delta^m_{B_n} $  as follows. 
 \begin{theorem}
\label{theorBD}
 The number of   $(i-1)$-dimensional faces of  $ \Delta^m_{B_n} $ is equal to 
  \[   \left( \begin{array}{c} mn+i \\ i \end{array} \right) \left( \begin{array}{c} n \\ i \end{array} \right) , \;\; i =0,\ldots, n.\]
    \end{theorem}   
   \vspace{0.1in}
   Given this, we  compute for $W=A_{n-1},B_n$  the $h$-vector of $\Delta^m_W $ and show that its $i$-th entry is the 
         $i$-th generalized Narayana number $ N^m_W(i) $   \cite{Atha3}. 
     These numbers sum up to the $m$-th generalized  Catalan number 
      associated to $ W $, defined in \cite{Atha2},  and admit several combinatorial interpretations. For instance they count 
          the regions  of the  arrangement  of hyperplanes $ (\alpha,x)=k, \; \alpha  \in \Phi^+,\;  k=0,\ldots,m $  
            in the fundamental chamber  that have $r-i$ walls of the form $(\alpha,x)=m $ separating them from 
      the fundamental alcove, the co-filtered chains of ideals of $\Phi^+ $ of length $m$  
    with $r-i$ indecomposable elements of rank $m$, as well as the 
     orbits of rank  $r-i$    of the action of $W$ on a certain quotient of the coroot lattice, where $r$ is the rank of the Weyl group $W$ and 
  $ \Phi^+$ is the set of    positive roots of $\Phi$. See \cite{Atha3} for further details. 
 
    \begin{corollary}
  \label{Narayana numbers}
   The i-th entry  of the $h$-vector of $ \Delta^m_W $ is equal to 
  \[  \frac {1}{i+1} \left( \begin{array}{c} n-1 \\i \end{array} \right)  \left( \begin{array}{c} mn \\i\end{array} \right) ,\;\; i=0,\ldots, n-1, \]
     \[ \left( \begin{array}{c} n\\i \end{array} \right) \left( \begin{array}{c} mn\\ i \end{array} \right) ,  \;\; i=0,\ldots,n,  \]
      in the cases  $ W  = A_{n-1}$ and $B_n$  respectively. 
  \end{corollary}

  The above numbers are exactly the  numbers $ N^m_W (i) $
 for the classical types $ A_{n-1}$ and $B_n $,  as  computed  in \cite[Section 5]{Atha3}. 
 We should point out that in the special case $W=B_n $ and $m=1$, Theorem \ref{theorBD} and Corollary \ref{Narayana numbers} 
  were first proved by Simion \cite[Section 3]{Si}. 

     For $m=1$ and $W=A_{n-1}$ or $B_n$,   $ \Delta^m_W$ is polytopal, i.e.  it is  isomorphic to the boundary complex of a polytope,  and therefore  it is shellable. 
  It follows from Corollary \ref{Narayana numbers} that for  $ m> 1 $ the $h$-vector of $ \Delta^m_W  $ is not symmetric, therefore $\Delta^m_W$ is not polytopal. 
   However,  we have the following theorem. 
   \begin{theorem}
  \label{shellability}
      The  complex $ \Delta^m_W$ is  shellable  for  the classical reflection groups $ W=A_{n-1},B_n $ and   $m \geq 1$. 
     \end{theorem}  
 \begin{corollary}
\label{geom description}
 The simplicial complex $ \Delta^m_W $  has the homotopy type of a  wedge of
        $ N^m_W (r) $ copies of the $(r-1)$-dimensional sphere, where $r$ is the rank of $W$,  for $ W= A_{n-1}, B_n$   respectively.    
   \end{corollary}
 
    Observe that $ \Delta^m_W$ is a flag complex (see \cite[Section III.4]{St}) for $ W=A_{n-1} $ and  $ B_n $. 
  Thus these complexes are examples of Cohen-Macaulay flag complexes.  
 Since the $h$-vector of a pure shellable complex forms an $M$-sequence (see \cite[Section II.3]{St})  we deduce the following corollary. 
\begin{corollary}
\label{M-vector}
  For the  reflection groups $ W=A_{n-1},B_n$  and   $m \geq 1$ the numbers 
   $ N^m_W (i) $ form an $M$-sequence.  
 \end{corollary}

 The statement of the previous corollary was conjectured by C. Athanasiadis (private communication) for all Weyl groups. 
 \vspace{0.1in}

 After this work was completed, it came to my attention that S. Fomin and N. Reading  \cite{FR} 
  have defined  and studied generalized cluster complexes for all finite root systems. 
   Our construction for types $A$ and $B$  agrees with their definition.

\vspace{0.1 in}
\emph{Aknowledgements}.  The definition of $\Delta^m_W$ in this paper was given by Victor Reiner in the case $W=A_{n-1}$ and Christos Athanasiadis 
   in the case $W=B_n$. I am grateful to both   for suggesting the study of the complexes $\Delta^m_W $ 
   and for indicating the possible connection with the generalized Narayana numbers. I also thank the anonymous referees 
  for their thorough comments. 


\section{Preliminaries}
\label{pre}
    In this section we introduce terminology and some elementary facts that  we will need later.   

   \vspace{0.1 in}
  \noindent
    {\bf $f$-vectors and $h$-vectors.}
   Let $E$ be a finite set. An { \em abstract simplicial complex}  on the ground set $E$ is a collection $ \Delta $ of subsets of $E$ such that $ F \subseteq  F' \in \Delta $ 
  implies $ F \in \Delta$. The set $ V = \{ v \in E : \{v\} \in \Delta \} $ is the set of {\em vertices}  of $ \Delta $.
   Every element  $ F $ of $\Delta $ is called a {\em face}. The  {\em dimension}  of $ F $ is one less than  its cardinality. 
    The {\em facets} of $ \Delta $ are the maximal faces with respect to inclusion   and $ \Delta $ is {\em  pure}  of  dimension  $d$  if every facet has
      dimension $ d $. We define the {\em induced subcomplex} of $\Delta$ on the vertex set $V' \subseteq V$ to be the simplicial complex whose faces are 
    the faces of $\Delta$ contained in $V'$. 

         Given a finite (abstract)  simplicial complex   $ \Delta$ of dimension $d-1$,
    let $ f(\Delta) = (f_{-1},f_0,\ldots,f_{d-1}) $ be the  {\em $f$-vector} of $\Delta$,  
    so that $ f_i $ is equal  to the number of  $i$-dimensional faces of $\Delta$. 
    The entry $ f_{-1} = 1 $ corresponds  to the empty face. 

   The $f$-vector determines the {\em  $h$-vector}  $h(\Delta) = (h_0,h_1,\ldots,h_d) $   through the  formula 
        \begin{equation}
     \label{relation}
       \sum_{i=0}^d f_{i-1} (x-1)^{d-i} = \sum_{i=0}^d  h_i x^{d-i}. 
     \end{equation}
     Equivalently we have 
    \begin{equation}
  \label{h-vector}
    h_k= \sum_{i=0}^k  (-1)^{k-i}  \left( \begin{array}{c} d-i \\d-k \end{array} \right) f_{i-1}
   \end{equation}
    \noindent  
    and 
    \begin{equation}
     \label{f-vector}
      f_{k-1}= \sum_{i=0}^k    h_i \left( \begin{array}{c} d-i \\k-i  \end{array} \right). 
   \end{equation}
\noindent
In particular, from (\ref{h-vector})  we have   
     \begin{equation}
   \label{EulerChar}
        (-1)^{d-1} h_d = -1+f_0 -f_1 + \cdots+ (-1)^{d-1} f_{d-1},  
  \end{equation}
     which is the reduced Euler charecteristic of $ \Delta$. 
       
  \vspace{0.1 in}
 \noindent 
 {\bf Vertex decomposability.}
   Let $\Delta $ be an abstract simplicial complex over the ground set $E$. 
 The  {\em cone}  of $ \Delta $ over a new vertex $ v $, denoted by $ \Delta * v $,  is the simplicial complex 
     on the ground set $ E \cup \{ v \} $  with facets $ F \cup \{ v \} $,  where $ F $ is a facet  of $ \Delta$. More generally, the 
  {\em simplicial join} of two  abstract simplicial complexes $ \Delta $ and $ \Delta'$  with disjoint ground sets $ E $ and $ E'$   respectively, 
     denoted by $ \Delta * \Delta' $, is the abstract simplicial complex on the ground set $ E \cup E' $ with facets $ F \cup F' $,  
     where $ F$ is a facet of $ \Delta$ and $F'$ is a facet of $\Delta'$. 
   Notice that the dimension of $ \Delta*\Delta'$ is one more than  the sum of the dimensions of $\Delta$ and $\Delta'$. 
    For $ A \subseteq E$  the {\em  deletion} of A from $ \Delta $ is  the set  $ \Delta \setminus A := \{ B \in \Delta : A \cap B = \emptyset \} $ 
    and the {\em link} of $ A $ in $\Delta$ 
is the set $ \Delta / A := \{ B \in \Delta : A \cap B = \emptyset , A \cup B \in \Delta \} $. Observe that these two operations commute. 

\begin{definition} { \rm(Provan and Billera \cite{Bi,Pro})}
\label{def vd}
 A simplicial complex  $\Delta $ is vertex decomposable if it is pure and it is either empty  or it has a vertex $v $ such that $\Delta \setminus  v $ and $\Delta / v $ are 
  vertex decomposable. 
\end{definition}
 
  Notice that if $\Delta$ is pure of dimension $k$ then $ \Delta / v $ is pure of dimension $k-1$. 
  If $ \Delta \setminus v $ is also pure, then either $ \dim( \Delta \setminus v) = k $ or $\Delta $ is a cone over $ \Delta \setminus v$, 
   where $ \Delta \setminus v = \Delta / v $ has dimension $k-1$.  

  It was proved in \cite{Bi,Pro} that  every vertex decomposable simplicial complex is shellable. 
  In our proofs we will use the following easy consequence of  Definition \ref{def vd}. 
\begin{lemma} \cite[Lemma 2.2]{Atha1}
\label{vd}
 Let $\Delta $ be a simplicial complex of dimension  $ d $  and $ v_1,v_2,\ldots,v_t $ be distinct vertices of $ \Delta $. If 
  {\rm  (i)}  $\Delta \setminus \{v_1,\ldots,v_t \} $   is vertex decomposable of dimension $d$  and 
   {\rm (ii)}  $ (\Delta /v_i) \setminus \{v_1,\ldots,v_{i-1}  \} $ is vertex decomposable of dimension $ d-1$
   for each $ 1 \leq i \leq t $ then $ \Delta $ is vertex decomposable. 
 \end{lemma}
  
 Induction on the number of vertices gives  the next lemma. 
 \begin{lemma} 
\label{cone}
  The simplicial join of vertex decomposable complexes is vertex decomposable. 
In particular, the cone over  a vertex decomposable complex is  vertex decomposable. 
\end{lemma}


\section{Enumerative results}
  \label{enumeration}

      We will   first  prove Theorem \ref{theorBD}  bijectively. The bijection is analogous to the one exhibited in  \cite[Theorem 1]{Pr}.    
        Let $ \Delta^m_{B_n}[i]$ denote the set of dissections in  $ \Delta^m_{B_n} $ with exactly $i$ $B$-diagonals, 
          $ 0 \leq i \leq n $. Unless otherwise stated, the vertices of  every $(2mn+2)$-gon will be labeled by 
       $1,$ $2,\ldots,$ $mn+1,$ ${\bar 1},$ ${\bar 2},\ldots,$ $\overline{mn+1} $ in anticlockwise order.  

    \begin{theorem}
      \label{typeB}
     There is a bijection $ \Phi_{m,n} $ between the set $ \Delta^m_{B_n} [i] $   and the cartesian product of the set 
      of all sequences $\{  1\leq a_1 \leq a_2 \leq \cdots \leq a_i \leq mn+1\} $ with  the set 
       of  $(\epsilon _1,\ldots,\epsilon _n) \in \{0,1\}^n $ with exactly $ i $ entries $ \epsilon_j $  equal to 1. 
       \end{theorem}

  \begin{corollary}
 The number of $ (i-1)$-dimensional faces of $ \Delta^m_{B_n}$ is equal to 
    \[   \left( \begin{array}{c} mn+i \\ i \end{array} \right) \left( \begin{array}{c} n \\ i \end{array} \right) .\]
     \end{corollary}
    
\vspace{0.1in}
\noindent
       {\em Proof of Theorem \ref{typeB}.}  
Given a $(2mn+2)$-gon $P$   we  define the {\em initial point} of each diagonal   as follows. 
    The  initial point of a diameter is  the one with the positive label and every nondiameter diagonal is 
      oriented in such a way that if you travel along it in this orientation  you have  the  center of the polygon on the left.  
    If $D$ is  a nondiameter $B$-diagonal then it has two segments with
   initial points $a$ and ${\bar a}$ respectively. We  define the initial point of $D$ to be the one with the positive label. 

  Let us abbreviate  $ \Delta_n = \Delta^m_{B_n} $ and $ \Delta_n[i] = \Delta^m_{B_n}[i] $. 
     We construct a bijection  $ \Phi_{m,n} $  by induction on $n$.  
        If $n=1 $ then $ \Delta_n $  consists of $ m+1 $ diameters, or equivalently $m+1$ isolated vertices. 
       For $i=0$ the map $\Phi_{m,1}$ sends the trivial dissection  to the empty set and 
           for  $i=1$   sends  each element $a_k {\bar a}_k $ of $\Delta_n[1]$  to the 
        pair $((a_k),(\epsilon_1))$ with $\epsilon_1=1$.  

    For the inductive step assume we have already constructed a bijection
          \[ \Phi_{m,n} :   \Delta_n[i]    \longrightarrow  \{1\leq a_1 \leq a_2 \leq \cdots \leq a_i \leq mn+1\} \]
    \[ \times \;   \{ (\epsilon _1,\epsilon _2,\ldots,\epsilon _n) \in \{0,1\}^n  \mbox{: exactly  $ i $  entries $ \epsilon_j $ are  equal to 1} \} \] 
     for an integer  $ m $ and  $ 0 \leq i \leq n$.  
    To define  the bijection $ \Phi_{m,n+1} $ let $D \in  \Delta_{n+1}[i] $.  The result is trivial for $i=0$. 
  If $i \geq 1$  let $ 1\leq a_1 \leq a_2 \leq \cdots \leq a_i  \leq m(n+1)+1 $ 
     be the initial points of the $B$-diagonals of $ D $ and $ \mathcal{I}=\{ a_1,\ldots,a_i, \bar{a}_1,\ldots,\bar{a}_i\}$.   
       We consider the action of $ \ZZ $ on the set of vertices of $P$ which is defined as follows. 
   If $a$ is a vertex of $P$ and $ k $ a positive integer then $a+k$ is the vertex of $P $ which is $k$ vertices apart,  
     in the anticlockwise order,  from $a$. 
  Consider the minimum  $a_j$  such that none of the vertices $a_j+1,\ldots,a_j+m $
      belongs to the set $\mathcal{I} $. 
    We set $ \epsilon_1 = 1 $ if $ \{a_j  (a_j +m+1)$, \;$ {\bar a}_j  ({\bar a}_j +m+1) \}$   is a $B$-diagonal of $D$
    and $ \epsilon _1 = 0 $ otherwise. Let  $P'$  be the  
      $ (2mn+2)$-gon  with vertices $ 1,\ldots,$$a_j,$$a_j+m+1,\ldots,$$m(n+1)+1,$ $\bar{1},\ldots,$ $\bar{a}_j,$ $\bar{a_j}+m+1,$ $ \ldots,$ $\overline{m(n+1)+1} $ 
    which is  dissected by the remaining $ i-\epsilon _1 $   $B$-diagonals of $D$ 
   and let $D'$ be the dissection in $P'$ defined by these diagonals. 
 Note that the initial point of each $B$-diagonal is the same, whether we consider $D'$  in $P$ or $P'$. 
  Therefore, the  initial points of the $B$-diagonals of $D'$ in $P'$  are $I_0= (a_1,\ldots,a_i) $ if $\epsilon_1=0$ 
  or $ I_1= (a_1,\ldots,a_{j-1},a_{j+1},\ldots,a_i)$ if $\epsilon_1=1$. 
  In either case, by induction  there is a  sequence  $I_{\epsilon_1}$ of initial points   and a  
    sequence $(\epsilon _2,\ldots,\epsilon _{n+1}) $ with $ i-\epsilon_1 $    
        entries $ \epsilon_j $ equal to $1$   assigned to  $D'$  by $\Phi_{m,n}$. 
           We define $ \Phi_{m,n+1}(D) $ to be the pair $ ((a_1,\ldots,a_i ) , (\epsilon _1,\epsilon _2,\ldots,\epsilon _{n+1}) ) $. 

    In order to see that $ \Phi_{m,n+1} $ is a bijection we construct its inverse.  We proceed by induction, 
     the case $n=1$ being obvious. Assume we have constructed the inverse map $ \Phi^{-1}_{m,n} $. 
To define $\Phi^{-1}_{m,n+1}$ consider the  sequences
 $ 1 \leq  a_1 \leq  a_2 \leq \cdots$$ \leq  a_i$ $\leq m(n+1)+1 $ and  $ (\epsilon_1,\ldots,\epsilon_{n+1}) \in \{0,1\}^{n+1}  $ with $ i $
    entries $ \epsilon_j $  equal to 1. For $i=0$ the zero sequence bijects  to the trivial dissection. 
 If $i \geq 1$ consider the set  $ \mathcal{I}=\{ a_1,\ldots,a_i, \bar{a}_1,\ldots,\bar{a}_i\}$.   
     We will define  a dissection $D \in \Delta_{n+1}[i] $  as follows. 
      Find the minimum  $a_j$  such that  none of the vertices  $a_j+1,\ldots,a_j+m$ of $P$ belongs to $ \mathcal{I}$. 
      Then $\{a_j (a_j +m+1)$,  $  {\bar a}_j ({\bar a}_j +m+1)\}$ is a 
 $B$-diagonal of $ D $ if $ \epsilon_1 = 1$ and it is not otherwise. 
  Let $P'$ be the  $(2mn+2)$-gon with vertices $1,\ldots,$ $a_j,$ $a_j+m+1,\ldots,$ $m(n+1)+1,$
         ${\bar a}_1,\ldots,$ ${\bar a}_j,$ $\bar {a_j}+m+1,\ldots,$ $\overline{m(n+1)+1}$. 
        Consider the pair $ ((a_1,\ldots,a_i) , (\epsilon_2,\ldots,\epsilon_{n+1})) $ if $\epsilon_1=0$ 
     or $ ((a_1,\ldots,a_{j-1},a_{j+1},\ldots,a_i) , (\epsilon_2,\ldots,\epsilon_{n+1})) $ if $\epsilon_1=1$
         with exactly $i-\epsilon_1 $ entries $\epsilon_j $ equal to 1. By induction  there is 
       a unique dissection $D'$  in  $P'$ with $i-\epsilon_1 $ $B$-diagonals assigned to it
          by $\Phi^{-1}_{m,n}$. 
      Therefore,  we define
     $ \Phi^{-1}_{m,n+1} (  (a_1,\ldots,a_i) , (\epsilon_1,\ldots,\epsilon_{n+1}) ) $ to be the dissection 
         consisting of $ D'$  with the $B$-diagonal $\{ a_j (a_j+m+1)$, $ {\bar a}_j (\bar{a_j}+m+1)\}$
 added if $ \epsilon_1=1$ 
     and  to be   the dissection $D'$ in the original polygon  if $ \epsilon_1=0$. 
  We leave it to the reader to check that this map is indeed the inverse of $ \Phi_{m,n+1} $. 
     $ \qed$

      \begin{example} 
  {\rm      In Figure $1 $  we have a {\rm 26}-gon with $ m=2 $, $ n=6 $ and $ i=4 $ diagonals. 
   The bijection of Theorem \ref{typeB} maps  the dissection onto the pair
  $ (a_1,a_2,a_3,a_4 )=(6,11,11,12) $ and $ (\epsilon _1,\epsilon _2,\epsilon _3,\epsilon _4,\epsilon _5,\epsilon_6)=
   (1,1,0,1,0,1) $.  }
    \end{example}

 \begin{remark}
  \label{one-diameter}
 The number of dissections in $ \Delta^m_{B_n}[i] $ having  one diameter is equal to 
 \[ \left( \begin{array}{c} mn+i \\ i \end{array} \right)  \left( \begin{array}{c} n-1 \\ i-1 \end{array} \right) .\]
 \end{remark}

Indeed, it  follows from Theorem \ref{typeB}  and the observation that 
an $i$-element dissection $D$ has a diameter if and only if the pair 
  $ ((a_1,\ldots,a_i),(\epsilon_1,\ldots, \epsilon_n))$
   assigned to $D$  satisfies $\epsilon_n=1$. 
      Alternatively, it can be derived from Theorem \ref{note}.

\begin{figure}[h]
\centering{\includegraphics [scale = 0.6 ]{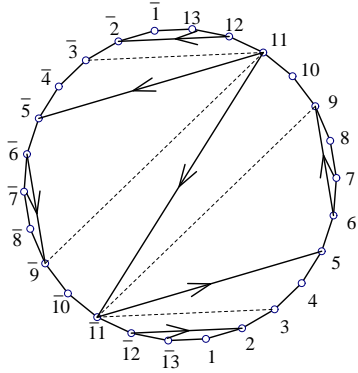}}
\caption{Example  3.3}
\end{figure}
 
   \noindent
  {\em {Proof of Corollary \ref{Narayana numbers}: }}
     It follows from formula (\ref{h-vector}), Theorems \ref{note} and \ref{theorBD} 
     and elementary calculations.      $\qed$

  \vspace{0.1in}
   From formula  (\ref{EulerChar}) we deduce  the following corollary.  
  \begin{corollary}
  For the classical reflection groups $W=A_{n-1},B_n$, the reduced Euler characteristic of $ \Delta^m_W $  equals  $ (-1)^{r-1} N^m_W (r) $, where $r$  is the rank of $W $. 
  \end{corollary}

    \section{Shellability}
In this section we prove Theorem \ref{shellability} by showing that $ \Delta^m_W $ is vertex decomposable for the classical reflection groups $W=A_{n-1},B_n$.

\begin{proposition}
\label{vd typeA}
 The simplicial complex $ \Delta^m_{A_{n-1}} $ is vertex decomposable.    
\end{proposition}
   
\begin{proof}
  We will prove a stronger statement. 
  Let $ P $ be the  $ (mn+2)$-gon and $ v $ a minimal diagonal, that is  a diagonal  which dissects $ P $ into an $(m+2)$-gon and an $(m(n-1)+2)$-gon.  
   Let  $a_1,a_2,\ldots,a_m $ be  the vertices  of the $ (m+2)$-gon in the anticlockwise order that are not endpoints of $v $, as in Figure $2$.  
     Denote by $\Delta^m_{A_{n-1}}(k) $ the complex  obtained from  $\Delta^m_{A_{n-1}} $ by deleting all 
   diagonals incident  to  $a_1,\ldots,a_k $ 
   and  let $ \Delta^m_{A_{n-1}}(0) = \Delta^m_{A_{n-1}} $.  
   We will prove that $\Delta^m_{A_{n-1}}(k) $ is vertex decomposable  and $(n-2)$-dimensional  for every 
     $ n $ and  $k, \;\; 0 \leq k \leq m $. For $k=0 $ this reduces to our proposition. 

  The statement is true for $n=1$ since the only element of  $ \Delta^m_{A_{n-1}} $ in this case is 
   the trivial  dissection. 
  We assume that $ \Delta^m_{A_l}(m-k) $ is vertex decomposable and $(l-1)$-dimensional for every $ 1\leq l \leq n-1$ and $ 0\leq k \leq m$ and 
   prove that  $ \Delta^m_{A_n}(m-k) $ is vertex decomposable and  $(n-1)$-dimensional  for every $k$, $ 0 \leq k \leq m $. We proceed by induction on $k$. 
       For $ k=0 $ the simplicial complex $ \Delta^m_{A_n}(m) $ is isomorphic to the cone $ v * \Delta^m_{A_{n-1}} (0) $, 
  which is vertex decomposable 
      and $(n-1)$-dimensional  by induction  and Lemma \ref{cone}. Assume that $ \Delta^m_{A_n}(m-k) $ is
   vertex decomposable of the same dimension. To  
      prove  the same for $ \Delta^m_{A_n}( m-k-1) $,   let $ v_1, \ldots,v_n $ be all the diagonals incident to  $ a_{m-k} $
    ordered by the clockwise ordering of their  endpoints other than $a_{m-k}$, as in Figure 2. 
     We  apply Lemma \ref{vd} on $ \Delta^m_{A_n}( m-k-1) $  and the set of diagonals  $ \{v_1,\ldots,v_n \} $.  
    Notice first that  $ \Delta^m_{A_n}(m-k-1) \setminus \{v_1,\ldots,v_n \} =$  $ \Delta^m_{A_n}(m-k) $, which   is 
   vertex decomposable and $(n-1)$-dimensional  by induction.  
   Therefore, condition (i) of Lemma \ref{vd} is verified. We then claim that 
   $  \Delta_i = $  $ \Delta^m_{A_n}(m-k-1)  / v_i \setminus \{v_1,\ldots,v_{i-1} \} $ is vertex decomposable  and $(n-2)$-dimensional  for every $ i=1,\ldots, n $. 
   Observe that the complex $ \Delta_i $ is obtained from 
     $ \Delta^m_{A_n}/ v_i $  by  deleting  the  diagonals $v_1,\ldots,v_{i-1} $ 
   and all  diagonals  incident to  $a_1,\ldots,a_{m-k-1} $. 
    Notice  that the diagonal $ v_i $ dissects   the $(2m(n+1)+2)$-gon 
  into an $ (mi+2)$-gon  on  which the points $ a_1,\ldots,a_{m-k}$ belong
and an $ (m(n+1-i)+2)$-gon. So, deleting  the preceding diagonals  from $\Delta^m_{A_n}/ v_i $ 
   is equivalent to deleting all  diagonals incident to  $a_1,\ldots,a_{m-k}$ in the $(mi+2)$-gon. 
  This implies that the  complex $ \Delta_i $ is the simplicial join  of two complexes isomorphic to 
$ \Delta^m_{A_{i-1}} (m-k) $ and $ \Delta^m_{A_{n-i}}( 0) $. This   is vertex decomposable  by Lemma \ref{cone} and induction 
  and has the desired dimension.  Therefore condition (ii) of Lemma \ref{vd} is satisfied and this completes the proof. 
\end{proof}

\begin{figure}[h]
\label{typeA}
\centering{\includegraphics [scale = 0.6 ]{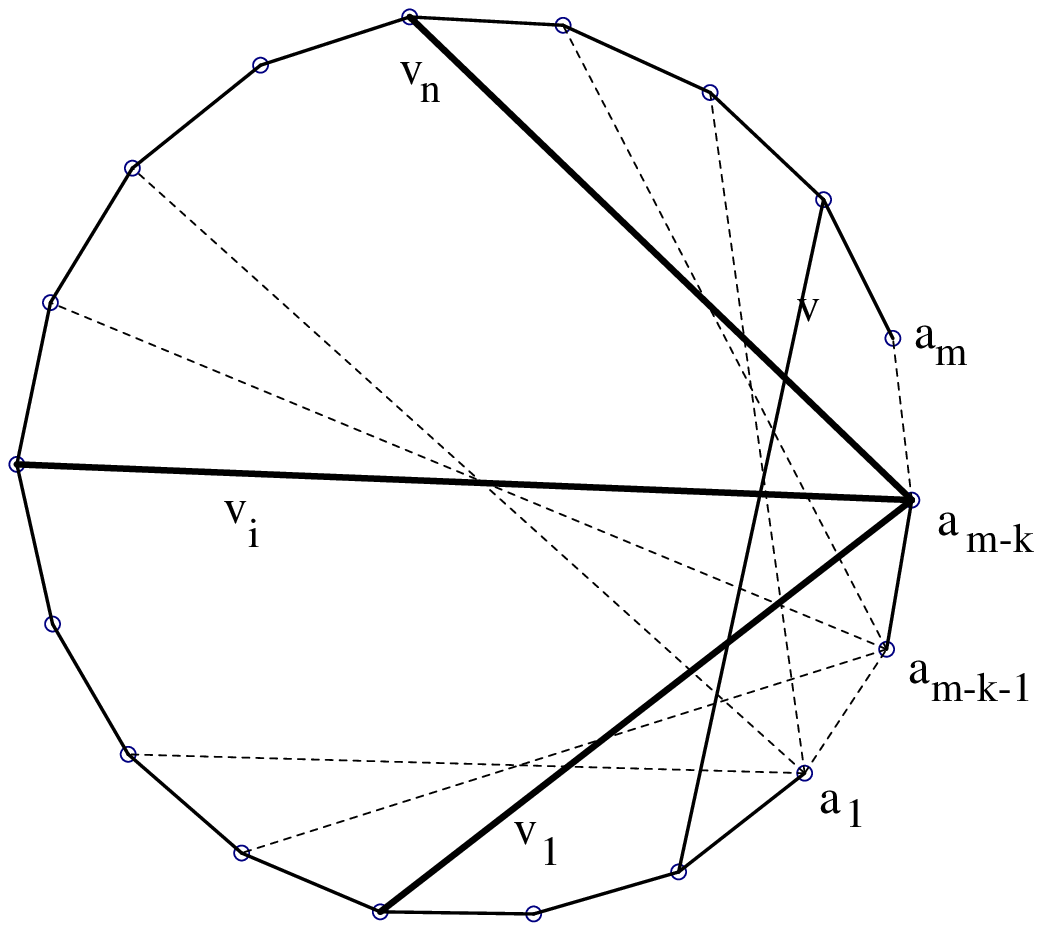}}
\caption{}
\end{figure}

    \begin{proposition}
\label{vdB}
   The simplicial complex $\Delta^m_{B_n} $ is vertex decomposable.  
   \end{proposition}

\begin{proof}
We will again prove the analogous stronger statement. 
  Let $ P $ be  the   $ (2mn+2)$-gon and $ v $ a minimal  $B$-diagonal, that is  a diagonal  which dissects $ P $ into
  a pair of  $(m+2)$-gons  and a centrally symmetric  $(2m(n-1)+2)$-gon.  
   Let  $a_1,a_2,\ldots,a_m $ be  the vertices in anticlockwise order of one of the $ (m+2)$-gons  that are not endpoints of $v $, as in Figure $3$. 
   Denote by $ \Delta^m_n(k) $ the complex obtained from  $ \Delta^m_{B_n} $   by  deleting
 all  diagonals incident to $a_1,\ldots,a_k $ 
   and  let  $ \Delta^m_n (0) = \Delta^m_{B_n} $. We will prove that $ \Delta^m_n(k) $ is vertex decomposable and $(n-1)$-dimensional for every 
$n$ and  $k,$ $ 0\leq k \leq m$. 

 The statement is true for $n=0$, since the only element of $ \Delta^m_{B_n} $ in this case is  the trivial  dissection.  
  We assume that $ \Delta^m_{l}(m-k) $ is vertex decomposable and $(l-1)$-dimensional for every $ 1 \leq l \leq n$, $ 0 \leq k \leq m$ and 
   we  prove that  $\Delta^m_{n+1}(m-k)$ is  $n$-dimensional and vertex decomposable for  $0 \leq k \leq m$. 
   We proceed by induction on $k$.  
   For $ k=0 $,   $ \Delta^m_{n+1}(m) $ is isomorphic to the cone of  $ \Delta^m_n(0) $ over $v$, 
 thus vertex decomposable and $n$-dimensional by induction and Lemma \ref{cone}.
   Next we assume   that  $ \Delta^m_{n+1}(m-k) $ is vertex decomposable and $n$-dimensional  and  prove the same for $ \Delta^m_{n+1}(m-k-1) $.  
 Let $ v_1,\ldots,v_{2n+1} $  be the $B$-diagonals incident  to  $ a_{m-k} $ 
ordered by the clockwise ordering of their  endpoints other than $a_{m-k}$. 
We apply Lemma \ref{vd} on $\Delta^m_{n+1}(m-k-1) $ and the set of these diagonals. 
First,   $ \Delta^m_{n+1}(m-k-1) \setminus \{ v_1,\ldots,v_{2n+1} \} = $  $  \Delta^m_{n+1} (m-k) $ is $n$-dimensional and  vertex decomposable by induction.  
   Thus, condition  (i) of Lemma \ref{vd}  holds. We then claim that $\Delta_i=$ $\Delta^m_{n+1}(m-k-1) / v_i $ $\setminus \{ v_1, \ldots, v_{i-1} \} $ is 
  $(n-1)$-dimensional and  vertex decomposable  for every $ 1 \leq i \leq 2n+1$. 
  Observe that $ \Delta_i $ is obtained from  $ \Delta^m_{B_{n+1}}/v_i $ by  deleting   the 
 $B$-diagonals $ v_1,\ldots,v_{i-1} $ and all  $B$-diagonals 
 incident to  $a_1,\ldots,a_{m-k-1}$. Denote by $\mathcal S$ the set of these $B$-diagonals. 
     We have to distinguish two cases. First, assume that $ i \leq n+1 $ (Figure 3(a)). 
    In this case $ v_i $ dissects the $ (2m(n+1)+2)$-gon into a pair of  $(mi+2)$-gons on which 
  the points $ a_1,\ldots,a_{m-k} $ lie
   and a centrally symmetric  $( 2m(n+1-i)+2)$-gon. 
  So deleting the $B$-diagonals of $\mathcal S$ 
from $\Delta^m_{B_{n+1}}/ v_i $ is equivalent to deleting
   all  $B$-diagonals incident to  $a_1,\ldots,a_{m-k}$ in the  $(mi+2)$-gons. 
Therefore  $\Delta_i $ is the simplicial join of  
$ \Delta^m_{n+1-i}(0) $ and   $ \Delta^m_{A_{i-1}}(i-1) $. 
    If $ i > n+1 $ (Figure 3(b)) then $v_i $ dissects the $(2m(n+1)+2)$-gon into a 
     centrally symmetric $ ( 2m(i-n-1)+2)$-gon and a pair of $ (( 2n+2-i)m+2)$-gons. 
  In this case  the points $a_1,\ldots,a_{m-k}$  belong to  the centrally symmetric polygon. 
 Therefore deleting the  $B$-diagonals of $\mathcal S$ from $\Delta^m_{B_{n+1}}/ v_i $
is equivalent to deleting all  $B$-diagonals incident to 
 $a_1,\ldots,a_{m-k}$ in the  $ ( 2m(i-n-1)+2)$-gon, forcing $ \Delta_i $ to be the simplicial join of 
     $ \Delta^m_{i-n-1}( 2(i-n-1)-1) $ and $ \Delta^m_{A_{2n+1-i}}(0)$.       
     In either cases, $\Delta_i $ is $(n-1)$-dimensional and vertex decomposable by induction,  Lemma \ref{vd} and Proposition \ref{vd typeA}. 
      Therefore, condition (ii) of Lemma \ref{vd}  is satisfied, which completes the induction. 
\end{proof}

\begin{figure}[h]
\label{third fig}
\centering{\includegraphics [scale = 0.7]{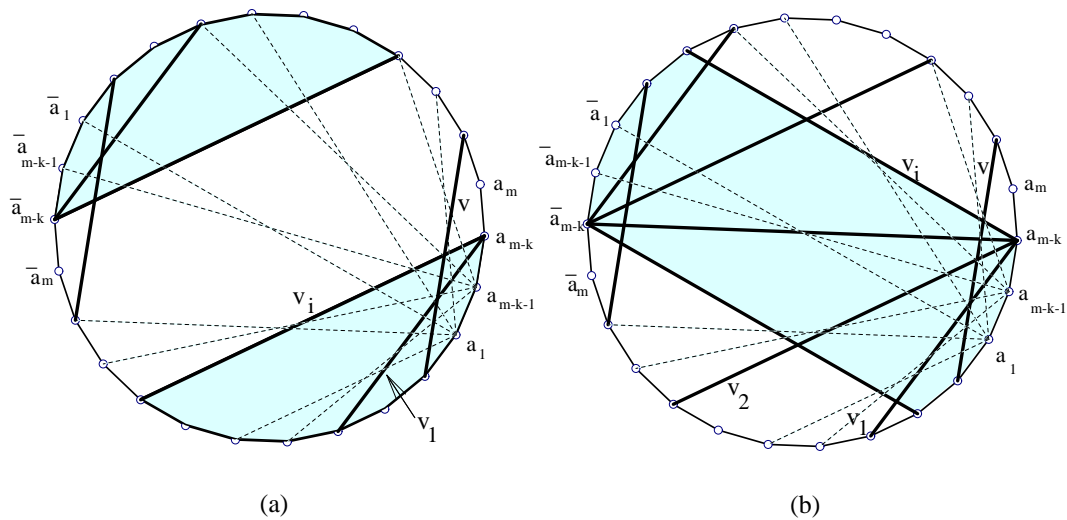}} 
\caption{The complex $\Delta_i$ when $ i \leq n $ and $ i>n $. 
Every diagonal represented by a dotted line comes in pair with a  centrally symmetric one which we have not drawn.}
\end{figure}


\end{document}